\documentclass[11pt]{amsart}
\usepackage{amssymb, amsthm, epsfig}
\usepackage[mathscr]{euscript}

\newtheorem{thm}{Theorem}[section]
\newtheorem{cor}[thm]{Corollary}
\newtheorem{lem}[thm]{Lemma}

\newtheorem{defnn}[thm]{Definition}
\newenvironment{definition}{\begin{defnn} \em}{\end{defnn}}
\newtheorem{remarkk}[thm]{Remark}

\newtheorem{examplee}[thm]{Example}

\newcommand{\scrg}{\mathcal{G}}
\newcommand{\scrc}{\mathcal{C}}

\newcommand{\scrf}{\mathcal{F}}

\newcommand{\bbz}{\mathbb{Z}}

\newcommand{\bbq}{\mathbb{Q}}
\newcommand{\bbn}{\mathbb{N}}

\newcommand{\G}{\Gamma}

\newcommand{\Ker}{\operatorname{Ker}}
\def\L{\Lambda}
\def\G{\Gamma}
\def\e{\epsilon}
\def\int{\mathrm{int}}

\title[Non-concordant knots having the same Seifert form]
{Infinite family of non-concordant knots\\ having the same Seifert
form}
\author{Taehee Kim}
\date{\today}
\address{Math Department--MS 136, Rice University,
6100 S.~Main St.~, Houston, TX 77005-1892, USA}
\email{tkim@rice.edu}
\def\subjclassname{\textup{2000} Mathematics Subject Classification}
\expandafter\let\csname
subjclassname@1991\endcsname=\subjclassname
\expandafter\let\csname
subjclassname@2000\endcsname=\subjclassname

\begin{document}
\begin{abstract}
By a recent result of Livingston, it is known that if a knot has a
prime power branched cyclic cover that is not a homology sphere,
then there is an infinite family of non-concordant knots having
the same Seifert form as the knot. In this paper, we extend this
result to the full extent. We show that if the knot has nontrivial
Alexander polynomial, then there exists an infinite family of
non-concordant knots having the same Seifert form as the knot. As
a corollary, no nontrivial Alexander polynomial determines a
unique knot concordance class. We use Cochran-Orr-Teichner's
recent result on the knot concordance group and Cheeger-Gromov's
von Neumann rho invariants with their universal bound for a
3-manifold.
\end{abstract}

\maketitle \vspace*{-2em}

\section{Introduction}
We work in the topologically locally flat category. A \emph{knot}
is an embedding of a circle into the 3-sphere. A knot is called
\emph{slice} if it bounds a (locally flat) 2-disk in the 4-ball.
For two knots $K_1$ and $K_2$, $K_1$ is said to be
\emph{concordant} to $K_2$ if $K_1\#-K_2$ is slice. Here the
symbol $\#$ denotes the connected sum operation and $-K$ denotes
the mirror image of $K$ with reversed orientation. This is an
equivalence relation. The equivalence classes (which are called
\emph{the concordance classes}) form an abelian group under the
connected sum operation. The group is called \emph{the (classical)
knot concordance group} and denoted by $\scrc$. In $\scrc$, the
identity is the class of slice knots. Levine \cite{l} constructed
an epimorphism $\phi : \scrc \rightarrow \scrg$ where $\scrg$
denotes the algebraic concordance group of Seifert forms modulo a
certain equivalence relation. The homomorphism $\phi$ maps the
concordance class represented by a knot to the algebraic
concordance class represented by Seifert forms of the knot. Jiang
\cite{j} showed that the kernel of $\phi$ is infinitely generated.
This implies that for each algebraic concordance class there are
infinitely many (mutually) non-concordant knots whose Seifert
forms represent that algebraic concordance class. But each
algebraic concordance class is also represented by infinitely many
distinct Seifert forms, and a question arises whether or not for a
given Seifert form there are non-concordant knots having
\emph{that} Seifert form. In fact, Jiang's examples have distinct
Seifert forms, hence his result does not give an answer to this
question. Recently Livingston \cite{li} made progress and gave a
partial answer under a condition on the Alexander polynomials.




\newtheorem*{liv1}{Theorem}
\begin{liv1}\cite[Theorem1.1]{li} If a knot $K$ has Seifert
form $V_K$ and its Alexander polynomial $\Delta_K(t)$ has an
irreducible factor that is not a cyclotomic polynomial $\phi_n$
with $n$ divisible by three distinct primes, then there is an
infinite family $\{K_i\}$ of non-concordant knots such that each
$K_i$ has Seifert form $V_K$.
\end{liv1}

In the above theorem the technical condition on the Alexander
polynomial is necessary since the theorem was proven by using
Casson-Gordon invariants. (For Casson-Gordon invariants, refer to
\cite{cg}.) More precisely, Casson-Gordon invariants are defined
via characters on the first homology of prime power branched
cyclic covers of knots and if every prime power branched cyclic
cover of the knot has the trivial first homology then all
Casson-Gordon invariants vanish. The following theorem due to
Livingston shows that a knot has a prime power branched cyclic
cover with nontrivial first homology under the given condition on
the Alexander polynomial. In the theorem, $\Delta_K(t)$ denotes
the Alexander polynomial of a knot $K$.

\newtheorem*{liv2}{Theorem}
\begin{liv2}\cite[Theorem1.2]{li} All prime power branched
cyclic covers of a knot $K$ are homology spheres if and only if
all nontrivial irreducible factors of $\Delta_K(t)$ are cyclotomic
polynomials $\phi_n(t)$ with $n$ divisible by three distinct
primes. All finite branched cyclic covers of $K$ are homology
spheres if and only if $\Delta_K(t) = 1$.
\end{liv2}

\noindent In addition to these results the author \cite{k} proved
that for each $n$ divisible by three distinct primes there exist
infinitely many non-concordant knots $K_i$ with $\Delta_{K_i}(t) =
(\phi_n(t))^2$ which have the same Seifert form. (In fact, in
\cite{k} the author showed that the knots $K_i$ are linearly
independent in the knot concordance group.)

In this paper we extend the above results to the full extent. The
main theorem is as follows.

\begin{thm}[Main Theorem]
\label{thm:main} If a knot $K$ has Seifert form $V_K$ and its
Alexander polynomial is not 1, then there is an infinite family
$\{K_i\}$ of non-concordant knots such that each $K_i$ has Seifert
form $V_K$.
\end{thm}

\noindent We note that by Freedman's work if $\Delta_K(t) = 1$
then $K$ is topologically slice \cite{f,fq}. (That is, the
concordance class of $K$ is the identity in $\scrc$.) On the other
hand, the main theorem implies that if $\Delta_K(t)$ is not 1 then
there are infinitely many non-concordant knots having the
Alexander polynomial $\Delta_K(t)$. Thus we have the following
corollary.

\begin{cor}
\label{cor:main} No nontrivial Alexander polynomial determines a
unique concordance class in the knot concordance group.
\end{cor}

In the proof of the main theorem we construct the knots $K_i$ by
performing \emph{satellite construction} on $K$. (This
construction is also called \emph{genetic modification} in
\cite{cot2}.) This construction is briefly reviewed in the next
section. To show that the $K_i$ are mutually non-concordant we use
Cochran-Orr-Teichner's filtration of the knot concordance group in
\cite{cot1} and von Neumann $\rho$-invariants defined by Cheeger
and Gromov \cite{chg}, which were applied as knot concordance
invariants first by Cochran, Orr, and Teichner in \cite{cot1}. In
particular, we use the fact that there is a universal bound for
von Neumann $\rho$-invariants for a fixed 3-manifold. More
precisely, for a fixed 3-manifold $M$, there exists a constant
$c_M$ such that $\left|\rho_\G^{(2)}(M,\psi)\right| \le c_M$ for
every representation $\psi : \pi_1(M) \rightarrow \Gamma$ where
$\Gamma$ is an arbitrary group \cite[Theorem 3.1.1]{r}. We remark
that in \cite{ct} Cochran and Teichner used this fact to show that
Cochran-Orr-Teichner's filtration of the knot concordance group is
highly nontrivial, that is, $\scrf_n/\scrf_{n.5}$ is nontrivial
for all $n \ge 2$.
\\

\noindent {\bf Acknowledgement} The author would like to thank Jae
Choon Cha for helpful conversations.

\section{Preliminaries}
Throughout this paper, we use the following convention. Unless
mentioned otherwise, integer coefficients are to be understood for
homology groups. The zero surgery on a knot $K$ in $S^3$ is
denoted by $M_K$. We use the same notation for a simple closed
curve and the homology (or homotopy) class represented by the
curve. We denote $\bbq[t,t^{-1}]$, the Laurant polynomial ring
with rational coefficients, by $\L$.

In this section we briefly review the machinery that will be used
in the proof of the main theorem. In \cite{cot1}, Cochran, Orr,
and Teichner established a filtration of the knot concordance
group $\{\scrf_n\}_{n\in \frac12\bbn_0}$ indexed by half-integers
where $\scrf_n$ is the subgroup of \emph{$(n)$-solvable} knots.
The definition of $(n)$-solvable knots ($n\in \bbn_0$) is as
follows. Recall that for a group $G$, $G^{(n)}$ denotes \emph{the
$n^{\text th}$ derived group of $G$} which is defined as follows:
Let $G^{(0)}\equiv G$, and inductively $G^{(n)} \equiv
[G^{(n-1)},G^{(n-1)}]$.
\begin{definition}
\label{def:n-solvable knots} A knot $K$ is called {\em
$(n)$-solvable} if $M_K$ bounds a spin $4$-manifold $W$ such that
the inclusion map $M_K \to W$ induces an isomorphism on the first
homology and such that $W$ admits an $(n)$-Lagrangian with
$(n)$-duals. This means that the intersection form (and the
self-intersection form) on
$H_2\left(W;\bbz[\pi_1(W)/\pi_1(W)^{(n)}]\right)$ pairs the
$(n)$-Lagrangian and the $(n)$-duals non-singularly and that their
images together freely generate $H_2(W)$. The 4-manifold $W$ is
called an \emph{$(n)$-solution for $K$} and we say $K$ is
\emph{$(n)$-solvable via $W$}.
\end{definition}

\noindent Similarly, we define $(n.5)$-solvable knots for $n\in
\bbn_0$. (An $(n.5)$-solution $W$ is required to admit an
$(n+1)$-Lagrangian with $(n)$-duals.) For more details, refer to
\cite[Definition 8.5 and Definition 8.7]{cot1}.

Cochran-Orr-Teichner showed that every slice knot is
$(n)$-solvable for all $n$ \cite[Remark 1.3.1]{cot1}. They detect
$(n.5)$-solvable knots, $n \in \bbn_0$, using von Neumann
$\rho$-invariants as follows.

\begin{thm}\cite[Theorem 4.2]{cot1}
\label{thm:cot main} Suppose $\G$ is an $(n)$-solvable
poly-torsion-free-abelian group. Let $\phi : \pi_1(M_K)
\rightarrow \G$ be a homomorphism. If $K$ is $(n.5)$-solvable via
a 4-manifold $W$ over which the coefficient system $\phi$ extends,
then $\rho^{(2)}_\G(M_K,\phi) = 0$.
\end{thm}

\noindent We explain the terminologies in the theorem. A group $G$
is called \emph{$(n)$-solvable} if $G^{(n+1)} = 1$. A group $G$ is
defined to be \emph{poly-torsion-free-abelian} (henceforth PTFA)
if it admits a normal series $1 = G_0 \lhd G_1\lhd \ldots\lhd G_m
= G$ such that the factors $G_{i+1}/G_i$ are torsion-free abelian.
For the von-Neumann $\rho$-invariant $\rho^{(2)}_\G(M_K,\phi)$,
refer to \cite[Section 5]{cot1} and \cite[Section 5]{cot2}.

In fact, the target group $\G$ which we will use for the proof of
the main theorem is a quotient group $G/G^{(n)}_r$ where
$G_r^{(n)}$ is \emph{the $n^\text{th}$ rational derived group} of
$G$ defined by Harvey \cite{h} as follows. Let $G_r^{(0)} \equiv
G$. For $n\ge 1$, define $G_r^{(n)} \equiv
\left[G_r^{(n-1)},G_r^{(n-1)}\right] P_{n-1}$ where
$$
P_{n-1} = \left\{g\in G_r^{(n-1)} \left| \right. g^k \in
\left[G_r^{(n-1)},G_r^{(n-1)}\right] \mbox{ for some }k\in
\bbz-\{0\}\right\}.
$$
The quotient $G_{r}^{(i)}/G_{r}^{(i+1)}$ is isomorphic to
$\left(G_{r}^{(i)}\left/ \left[G_{r}^{(i)},G_{r}^{(i)}\right]
\right. \right) \left/ \left\{
\bbz\!-\!\operatorname{torsion}\right\} \right. $ for all $i\ge 0$
\cite[Lemma 3.5]{h}. Harvey showed the quotient $G/G_r^{(n+1)}$ is
PTFA \cite[Corollary 3.6]{h}, and one easily sees that
$G/G_r^{(n+1)}$ is $(n)$-solvable.

To construct the knots $K_i$ in the main theorem we use
\emph{satellite construction} (or \emph{genetic modification})
explained as follows. Let $K$ be a knot and $\eta$ be an unknot in
$S^3$ which is disjoint from $K$. Let $J$ be another knot. Take
the union of the exterior of $\eta$ in $S^3$ and the exterior of
$J$ in $S^3$ along the common boundary (which is homeomorphic to a
torus) such that a meridian of $\eta$ is identified with a
longitude of $J$ and a longitude of $\eta$ with a meridian of $J$.
The resulting ambient manifold is homeomorphic to $S^3$. The image
of $K$ under this construction is denoted by $K(\eta,J)$ and we
say $K(\eta,J)$ is obtained by performing satellite construction
on $K$ via $\eta$ and $J$. If we let $D$ be an embedded disk in
$S^3$ bounded by $\eta$, then this construction is equivalent to
tying all the strands of $K$ transversally passing through $D$
into $J$. For more details, refer to \cite{cot2}. This
construction can be generalized to the case when we have a trivial
link $\{\eta_1, \ldots, \eta_n\}$ which is disjoint from $K$ and
auxiliary knots $\{J_1, \ldots, J_n\}$ by iterating the above
process. In this case the resulting knot is denoted by
$K(\{\eta_1, \ldots, \eta_n\}, \{J_1, \ldots, J_n\})$.

\section{Proof of Theorem~\ref{thm:main}}
Let $F$ be a Seifert surface of a knot $K$ with $\Delta_K(t) \ne
1$ and $V_K$ an associated Seifert form. The Seifert surface $F$
can be thought of as a disk with $2g$ bands where $g$ is the genus
of $F$. Let $i\in\bbn$. Let $\eta_i^n$, $1\le n \le 2g$, be a
trivial link in $S^3$ which is disjoint from $F$ such that the
$n^{\text th}$ component $\eta_i^n$ links the $n^{\text th}$ band
of $F$ once and does not link the other bands. It is known that
$\eta_i^n$, $1\le n \le 2g$, generate the rational Alexander
module $H_1(M_K;\L)$. (For example, see \cite{ro}.)

By \cite[Theorem 3.1.1]{r}, there exists a constant $c$ such that
$\left|\rho_\G^{(2)}(M_{K\#-K},\phi)\right| \le c$ for every
representation $\phi : \pi_1(M_{K\#-K}) \rightarrow \G$ where $\G$
is an arbitrary group. Let $J_1$ be a knot with vanishing Arf
invariant such that $\rho_\bbz^{(2)}(J_1) > c$. Here
$\rho_\bbz^{(2)}(J_1)$ denotes the von Neumann $\rho$-invariant
$\rho_\bbz^{(2)}(M_{J_1}, \phi)$ where $\phi : \pi_1(M_{J_1})
\rightarrow \bbz$ is the abelianization. Note that
$\rho_\bbz^{(2)}(J_1)$ is equal to the integral of the
Levine-Tristram signatures of $J_1$, integrated over the circle
normalized to length one \cite[Proposition 5.1]{cot2}.
Inductively, we define $J_{i+1}$ to be a knot with vanishing Arf
invariant such that $\rho_\bbz^{(2)}(J_{i+1}) > c +
2g\cdot\rho_\bbz^{(2)}(J_i)$. These $J_i$ can be easily found by
taking the connected sum of suitably many even copies of a
left-handed trefoil. For each $i\in\bbn$, let $J_i^n$ be a copy of
$J_i$ for $1 \le n \le 2g$. That is, $J_i^n \equiv J_i$, $1 \le n
\le 2g$.

Now let $K_i \equiv K(\{\eta_i^1, \ldots, \eta_i^{2g}\}, \{J_i^1,
\ldots, J_i^{2g}\})$, the knot resulting from satellite
construction. Since $\eta_i^n$, $1\le n \le 2g$, lie in the
complement of $F$ in $S^3$, $K_i$ have the same Seifert form $V_K$
as $K$. We prove $K_i$ are mutually non-concordant.

Fix $i$ and $j$ ($i< j$), and suppose that $K_i$ and $K_j$ are
concordant. That is, $K_i\# -K_j$ is slice. Observe that
\begin{equation*}
\label{eqn:genetic} K_i\#-K_j = (K\#-K)(\{\eta_i^1, \ldots,
\eta_i^{2g}, \eta_j^1, \ldots, \eta_j^{2g}\}, \{J_i^1, \ldots,
J_i^{2g}, -J_j^1, \ldots, -J_j^{2g}\}).
\end{equation*}

\noindent Here we abuse the notation so that $\eta^n_j$, $1\le
n\le 2g$, also denote the mirror images of $\eta^n_j$. Let $M
\equiv M_{K\#-K}$ and $M' \equiv M_{K_i\#-K_j}$.

We construct a cobordism $C$ between $M$ and $M'$ as follows.
Choose a $(0)$-solution $W_i$ (which is a spin 4-manifold) for
$J_i$. (Since $J_i$ has vanishing Arf invariant, it is
$(0)$-solvable. See \cite[Remark 1.3.2]{cot1}.) By doing surgery
along $\pi_1(W_i)^{(1)}$, we may assume tat $\pi_1(W_i) \cong
\bbz$. Similarly, we choose a $(0)$-solution $V_j$ for $-J_j$. Let
$W_i^n \equiv W_i$ and $V_j^n \equiv V_j$ for $1 \le n \le 2g$.
Take $M \times [0,1]$ and the disjoint union
$\left(\coprod_{n=1}^{2g}W_i^n\right) \coprod
\left(\coprod_{n=1}^{2g}V_j^n\right)$. To form C, for each $n$
identify the solid torus in $\partial W_i^n = (S^3\setminus
J_i^n)\cup S^1\times D^2$ with a tubular neighborhood of $\eta_i^n
\times \{1\}$ in $M\times \{1\}$ such that a meridian of $J_i^n$
is identified with a longitude of $\eta_i^n$ and a longitude of
$J_i^n$ with a meridian of $\eta_i^n$, and identify the solid
torus in $\partial V_j^n = (S^3 \setminus -J_j^n)\cup S^1 \times
D^2$ with a tubular neighborhood of $\eta_j^n \times \{1\}$ in
$M\times \{1\}$ similarly. One sees that $\partial_{-}C = M$ and
$\partial_{+}C = M'$. Moreover one sees that $C$ is spin.

Since $K_i\# -K_j$ is slice, $K_i\#-K_j$ is $(1.5)$-solvable by
\cite[Remark 1.3.1]{cot1}. Let $W'$ be a $(1.5)$-solution for
$K_i\#-K_j$. In particular, $\partial W' = M'$. Let $W$ be the
union of $C$ and $W'$ along their common boundary $M'$. Hence $W$
is a 4-manifold with $\partial W = M$.

\begin{lem}
\label{lem:1-solution} The 4-manifold $W$, which is constructed as
above, is a $(1)$-solution for $K\#-K$.
\end{lem}

The proof of the above lemma is postponed. Let $\G \equiv \pi_1(W)
/\pi_1(W)_r^{(2)}$. Note that $\G$ is a $(1)$-solvable PTFA group
by \cite[Corollary 3.6]{h}. Let $\phi : \pi_1(W) \rightarrow \G$
be the projection homomorphism. Note that $M'$, $M_{J^n_i}$,
$M_{-J^n_j}$, and $W'$ are subspaces of $W$, hence $\phi$ can be
restricted to the corresponding fundamental groups. By
\cite[Proposition 3.2]{cot2},
\begin{eqnarray*}
\rho_\G^{(2)}(M,\phi) & = & \rho_\G^{(2)}(M',\phi|_{\pi_1(M')})\\
& &{}+\sum_{n=1}^{2g}\rho_\G^{(2)}(M_{J^n_i},
\phi|_{\pi_1(M_{J^n_i})})\\
& &{}+\sum_{n=1}^{2g}\rho_\G^{(2)}(M_{-J^n_j},
\phi|_{\pi_1(M_{-J^n_j})}).
\end{eqnarray*}
In the above equation, $\rho_\G^{(2)}(M',\phi|_{\pi_1(M')}) = 0$
by Theorem~\ref{thm:cot main} since $\phi|_{\pi_1(M')}$ extends
over $(1.5)$-solution $W'$. Note that $\phi|_{\pi_1(M_{J^n_i})}$
factors through $\pi_1(W^n_i)$ which is isomorphic to $\bbz$ for
each $n$. If $\phi(\eta_i^n) = e$, the identity element in $\G$,
then $\rho_\G^{(2)}(M_{J^n_i}) = 0$. If $\phi(\eta_i^n) \ne e$,
then the image of $\phi|_{\pi_1(M_{J^n_i})}$ is isomorphic to
$\bbz$ and $\rho_\G^{(2)}(M_{J^n_i}) = \rho_\bbz^{(2)}(J^n_i)$,
which is defined in the previous section, by \cite[Proposition
5.13]{cot1}. We obtain similar results for
$\rho_\G^{(2)}(M_{-J^n_j})$. Now let $\e_i^n \equiv 0$ if
$\phi(\eta_i^n) = e$ and $\e_i^n \equiv 1$ otherwise, $1\le n \le
2g$. Define $\e_j^n$, $1\le n\le 2g$, similarly. Then we have the
following equation.

$$
\rho_\G^{(2)}(M,\phi) =
\sum_{n=1}^{2g}\e_i^n\rho_\bbz^{(2)}(J_i^n)
-\sum_{n=1}^{2g}\e_j^n\rho_\bbz^{(2)}(J_j^n).
$$

We claim that $\e_i^n \ne 0$ for some $n$ or $\e_j^n \ne 0$ for
some $n$. One sees that $\eta_i^n$ together with $\eta_j^n$, $1\le
n \le 2g$, generate the rational Alexander module $H_1(M;\L)$.
(This is obvious since $H_1(M;\L)$ is isomorphic to
$H_1(M_K;\L)\oplus H_1(M_{-K};\L)$.) Since $\Delta_K(t) \ne 1$,
$H_1(M;\L)$ is not trivial. Hence $K\#-K$ has the
\emph{nonsingular} rational Blachfield form $B\ell : H_1(M;\L)
\times H_1(M;\L) \rightarrow \bbq(t)/\L$. Let $i_* : H_1(M;\L)
\rightarrow H_1(W;\L)$ be the homomorphism induced by the
inclusion. Since $P \equiv \Ker(i_*)$ is self-annihilating by
\cite[Theorem 4.4]{cot1} (that is, $P = P^\bot$) and $B\ell$ is
nonsingular, $i_*$ is not a trivial homomorphism. Hence
$i_*(\eta_i^n) \ne 0$ for some $n$ or $i_*(\eta_j^n) \ne 0$ for
some $n$ in $H_1(W;\L)$. Since $W$ is a $(1)$-solution for
$K\#-K$, $H_1(W) \cong \bbz$. This implies that $\pi_1(W)^{(1)}_r
= \pi_1(W)^{(1)}$. Hence
$$
\pi_1(W)^{(1)}/\pi_1(W)^{(2)}_r \otimes_\bbz \bbq \cong
\pi_1(W)^{(1)}/\pi_1(W)^{(2)} \otimes_\bbz \bbq \cong H_1(W;\L).
$$
The first isomorphism holds by \cite[Lemma 3.5]{h}. Thus
$\phi(\eta_i^n) \ne e$ or $\phi(\eta_j^n) \ne e$ for some $n$ in
$\pi_1(W)^{(1)}/\pi_1(W)_r^{(2)}$ which is a subgroup of $\G$, and
this proves the claim.

Now suppose $\e_j^n \ne 0$ for some $n$. By our choice of $J_i$
and $J_j$,
$$
\rho_\G^{(2)}(M,\phi) \le 2g\cdot\rho_\bbz^{(2)}(J_i) -
\rho_\bbz^{(2)}(J_j) < -c,
$$
which is a contradiction. If $\e_j^n = 0$ for all $n$, then
$\e_i^n \ne 0$ for some $n$ by the above claim. Then
$$
\rho_\G^{(2)}(M,\phi) \ge \rho_\bbz^{(2)}(J_i) > c,
$$
which is also a contradiction. Therefore, to complete the proof we
only need to prove Lemma~\ref{lem:1-solution} and a proof is given
below.
\\ \\
\emph{Proof of Lemma~\ref{lem:1-solution}} : We follow a course of
the proof for a more general case in \cite{ct}. Using
Mayer-Vietoris sequence observe that
$$
H_1(M) \cong H_1(C) \cong H_1(M') \cong H_1(W') \cong H_1(W) \cong
\bbz.
$$
Again using Mayer-Vietoris sequence one sees that
$$
H_2(C) \cong \left(\oplus_{n=1}^{2g}H_2\left(W_i^n\right)\right)
\oplus \left(\oplus_{n=1}^{2g}H_2\left(V_j^n\right)\right) \oplus
H_2(M)
$$
and observe that $H_2(W) \cong \left(H_2(C) \oplus H_2(W')\right)
\left/\right.(p_*,q_*)(H_2(M'))$ where $p_*$ and $q_*$ are induced
by inclusions $p : M'\to C$ and $q : M'\to W'$, respectively.
Since $H^1(W') \rightarrow H^1(M')$ is an isomorphism, $H_3(W',M')
\rightarrow H_2(M')$ is an isomorphism by duality. Thus the
homomorphism $q_* : H_2(M') \rightarrow H_2(W')$ is the trivial
homomorphism. Observe that $H_2(M) \cong H_2(M') \cong \bbz$ and
they are generated by a capped-off Seifert surface of $K\#-K$ and
its image under satellite construction, respectively. Moreover
$p_* : H_2(M') \rightarrow H_2(C)$ maps the generator of $H_2(M')$
to the generator of $H_2(M)$. Hence
$$
H_2(W) \cong \left(\oplus^{2g}_{n=1}H_2\left(W_i^n\right)\right)
\oplus \left(\oplus^{2g}_{n=1}H_2\left(V_j^n\right)\right) \oplus
H_2(W').
$$
Observe that $\pi_1(W')^{(1)}$ maps into $\pi_1(W)^{(1)}$ by the
homomorphism induced by the inclusion. Also $\pi_1(W_i^n)$ and
$\pi_1(V_j^n)$ map into $\pi_1(W)^{(1)}$ by the homomorphisms
induced by the inclusions since $\eta_i^n$ and $\eta_j^n$ lie in
$\pi_1(W)^{(1)}$ and they generate $\pi_1(W_i^n)$ and
$\pi_1(V_j^n)$ (which are isomorphic to $\bbz$), respectively. Now
using naturality of equivariant intersection forms, one sees that
$(0)$-Lagrangians and $(0)$-duals for $W_i^n$ and $V_j^n$ and a
$(1)$-Lagrangian and $(1)$-duals for $W'$ consist of a
$(1)$-Lagrangian and $(1)$-duals for $W$. Finally, $W'$ has two
possible spin structures, and a spin structure on $W'$ can be
chosen such that $W$ is spin. This completes the proof.

\end{document}